\documentclass[aos,MSNbibl,nameyear,dvips]{arximspdf}
\usepackage{graphicx}
\doi{10.1214/13-AOS1123} 
\volume{41}
\issue{1}
\pubyear{2013}
\firstpage{1716}
\lastpage{1741}

\makeatletter

\newproclaim{definition}{Definition}

\newcommand{\rrvert}{\vert}
\newcommand{\llvert}{\vert}

\newtheorem{lemm}{Lemma}

\newcommand{\argmin}{\mathop{\arg\min}}

\newcommand{\cal}{\mathcal}

\newcommand{\xp}{\mathbf{x}_p}
\newcommand{\bps}{{\beta_{p*}}}
\newcommand{\bfX}{\mathbf{X}}
\newcommand{\bfI}{\mathbf{I}}
\newcommand{\bfy}{\mathbf{y}}
\newcommand{\bfXp}{\mathbf{X}_{-p}}
\newcommand{\con}{\mid}
\newcommand{\bF}{\mathbf{F}}
\newcommand{\bP}{\mathbf{P}}
\newcommand{\bft}{{\bolds\theta}}
\newcommand{\bS}{{\bolds\Sigma}}
\newcommand{\bfb}{{\bolds\beta}}

\makeatother

\begin{document}
\begin{frontmatter}

\title{Uniformly most powerful Bayesian tests}
\runtitle{Uniformly most powerful Bayesian tests}

\begin{aug}
\author[A]{\fnms{Valen E.} \snm{Johnson}\corref{}\thanksref{t1}\ead[label=e1]{vjohnson@stat.tamu.edu}}
\runauthor{V. E. Johnson}
\affiliation{Texas A\&M University}
\address[A]{Department of Statistics\\
Texas A\&M University\\
3143 TAMU\\
College Station, Texas 77843-3143\\
USA\\
\printead{e1}} 
\end{aug}

\thankstext{t1}{Supported by Award
Number R01 CA158113 from the National Cancer Institute.}

\received{\smonth{8} \syear{2012}}
\revised{\smonth{4} \syear{2013}}

%
\begin{abstract}
Uniformly most powerful tests are statistical hypothesis tests that
provide the greatest power against a fixed null hypothesis among all
tests of a given size. In this article, the notion of uniformly most
powerful tests is extended to the Bayesian setting by defining
uniformly most powerful Bayesian tests to be tests that maximize the
probability that the Bayes factor, in favor of the alternative
hypothesis, exceeds a specified threshold. Like their classical
counterpart, uniformly most powerful Bayesian tests are most easily
defined in one-parameter exponential family models, although extensions
outside of this class are possible. The connection between uniformly
most powerful tests and uniformly most powerful Bayesian tests can be
used to provide an approximate calibration between $p$-values and Bayes
factors. Finally, issues regarding the strong dependence of resulting
Bayes factors and $p$-values on sample size are discussed.
\end{abstract}

%

\end{frontmatter}

\section{Introduction}\label{sec1}
Uniformly most powerful tests (UMPTs) were proposed by Neyman and
Pearson in a series of articles published nearly a century ago [e.g.,
Neyman and Pearson (\citeyear{neyman28,neyman33}); see
\citet{lehmann05} for a comprehensive review of the subsequent
literature]. They are defined as statistical hypothesis tests that
provide the greatest power among all tests of a given size. The goal of
this article is to extend the classical notion of UMPTs to the Bayesian
paradigm through the definition of uniformly most powerful Bayesian
tests (UMPBTs) as tests that maximize the probability that the Bayes
factor against a fixed null hypothesis exceeds a specified threshold.
This extension is important from several perspectives.

From a classical perspective, the outcome of a hypothesis test is a
decision either to reject the null hypothesis or not to reject the null
hypothesis.
This approach to hypothesis testing is closely related to Popper's
theory of critical rationalism, in which scientific theories are never
accepted as being true, but instead are only subjected to increasingly
severe tests [e.g., \citet{popper59,mayo06}]. Many scientists and
philosophers, notably Bayesians, find this approach unsatisfactory for
at least two reasons [e.g., \citet{jeffreys39,howson05}]. First, a
decision not to reject the null hypothesis provides little quantitative
information regarding the truth of the null hypothesis. Second, the
rejection of a null hypothesis may occur even when evidence from the
data strongly support its validity. The following two examples---one
contrived and one less so---illustrate these concerns.

The first example involves a test for the distribution of a random
variable $X$ that can take values 1, 2 or 3; cf. \citet{berger88}. The
probability of each outcome under two competing statistical hypotheses
is provided in Table~\ref{bergertab}. From this table, it follows that
a most powerful test can be defined by rejecting the null hypothesis
when $X=2$ or 3. Both error probabilities of this test are equal to 0.01.

\begin{table}
\tablewidth=250pt
\caption{Probabilities of a random variable under competing
hypotheses}\label{tabl1}\label{bergertab}
\begin{tabular*}{\tablewidth}{@{\extracolsep{\fill}}lccc@{}}
\hline
$\bolds{X}$ & \textbf{1} & \textbf{2} & \textbf{3} \\
\hline
Null hypothesis & 0.99 & 0.008 & 0.001 \\
Alternative hypothesis & 0.01 & 0.001 & 0.989 \\
\hline
\end{tabular*}
\end{table}

Despite the test's favorable operating characteristics, the rejection
of the null hypothesis for $X=2$ seems misleading: $X=2$ is 8 times
more likely to be observed under the null hypothesis than it is under
the alternative. If both hypotheses were assigned equal odds a priori,
the null hypothesis is rejected at the 1\% level of significance even
though the posterior probability that it is true is 0.89. As discussed
further in Section~\ref{largesample}, such clashes between
significance tests and Bayesian posterior probabilities can occur in
variety of situations and can be particularly troubling in large sample
settings.

The second example represents a stylized version of an early phase
clinical trial. Suppose that a standard treatment for a disease is
known to be successful in 25\% of patients, and that an experimental
treatment is concocted by supplementing the standard treatment with the
addition of a new drug. If the supplemental agent has no effect on
efficacy, then the success rate of the experimental treatment is
assumed to remain equal to 25\% (the null hypothesis). A single arm
clinical trial is used to test this hypothesis. The trial is based on a
one-sided binomial test at the 5\% significance level. Thirty patients
are enrolled in the trial.

If $y$ denotes the number of patients who respond to the experimental
treatment, then the critical region for the test is $y \geq12$. To
examine the properties of this test, suppose first that $y=12$, so that
the null hypothesis is rejected at the 5\% level. In this case, the
\emph{minimum} likelihood ratio in favor of the null hypothesis is
obtained by setting the success rate under the alternative hypothesis
to $12/30 = 0.40$ (in which case the power of the test is 0.57). That
is, if the new treatment's success rate were defined a priori to be
0.4, then the likelihood ratio in favor of the null hypothesis would be
%
%
\begin{equation}
L_{\min} = \frac{0.25^{12} 0.75^{18}}{ 0.4^{12} 0.6^{18}} = 0.197.
\end{equation}
For any other alternative hypothesis, the likelihood ratio in favor of
the null hypothesis would be larger than 0.197 [e.g., \citet
{edwards63}]. If equal odds are assigned to the null and alternative
hypothesis, then the posterior probability of the null hypothesis is
\emph{at least} 16.5\%. In this case, the null hypothesis is rejected at
the 5\% level of significance even though the data support it. And, of
course, the posterior probability of the null hypothesis would be
substantially higher if one accounted for the fact that a vast majority
of early phase clinical trials fail.

Conversely, suppose now that the trial data provide clear support of
the null hypothesis, with only 7 successes observed during the trial.
In this case, the null hypothesis is not rejected at the 5\% level, but
this fact conveys little information regarding the relative support
that the null hypothesis received. If the alternative hypothesis
asserts, as before, that the success rate of the new treatment is 0.4,
then the likelihood ratio in favor of the null hypothesis is 6.31; that
is, the data favor the null hypothesis with approximately 6:1 odds. If
equal prior odds are assumed between the two hypotheses, then the
posterior probability of the null hypothesis is 0.863. Under the
assumption of clinical equipoise, the prior odds assigned to the two
hypotheses are assumed to be equal, which means the only controversial
aspect of reporting such odds is the specification of the alternative
hypothesis.

For frequentists, the most important aspect of the methodology reported
in this article may be that it provides a connection between
frequentist and Bayesian testing procedures. In one-parameter
exponential family models with monotone likelihood ratios, for example,
it is possible to define a UMPBT with the same rejection region as a
UMPT. This means that a Bayesian using a UMPBT and a frequentist
conducting a significance test will make identical decisions on the
basis of the observed data, which suggests that either interpretation
of the test may be invoked. That is, a decision to reject the null
hypothesis at a specified significance level occurs only when the Bayes
factor in favor of the alternative hypothesis exceeds a specified
evidence level. This fact provides a remedy to the two primary
deficiencies of classical significance tests---their inability to
quantify evidence in favor of the null hypothesis when the null
hypothesis is not rejected, and their tendency to exaggerate evidence
against the null when it is. Having determined the corresponding UMPBT,
Bayes factors can be used to provide a simple summary of the evidence
in favor of each hypothesis.

For Bayesians, UMPBTs represent a new objective Bayesian test, at least
when objective Bayesian methods are interpreted in the broad sense. As
\citet{berger06} notes, ``there is no unanimity as to the definition of
objective Bayesian analysis\ldots'' and ``many Bayesians object to the
label `objective Bayes,''' preferring other labels such as
``noninformative, reference, default, conventional and
nonsubjective.'' Within this context, UMPBTs provide a new form of
default, nonsubjective Bayesian tests in which the alternative
hypothesis is determined so as to maximize the probability that a Bayes
factor exceeds a specified threshold. This threshold can be specified
either by a default value---say 10 or 100---or, as indicated in the
preceding discussion, determined so as to produce a Bayesian test that
has the same rejection region as a classical UMPT. In the latter case,
UMPBTs provide an objective Bayesian testing procedure that can be used
to translate the results of classical significance tests into Bayes
factors and posterior model probabilities. By so doing, UMPBTs may
prove instrumental in convincing scientists that commonly-used levels
of statistical significance do not provide ``significant'' evidence
against rejected null hypotheses.

Subjective Bayesian methods have long provided scientists with a formal
mechanism for assessing the probability that a standard theory is true.
Unfortunately, subjective Bayesian testing procedures have not
been---and will likely never be---generally accepted by the scientific
community. In most testing problems, the range of scientific opinion
regarding the magnitude of violations from a standard theory is simply
too large to make the report of a single, subjective Bayes factor
worthwhile. Furthermore, scientific journals have demonstrated an
unwillingness to replace the report of a single $p$-value with a range
of subjectively determined Bayes factors or posterior model probabilities.
%

Given this reality, subjective Bayesians may find UMPBTs useful for
communicating the results of Bayesian tests to non-Bayesians, even when
a UMPBT is only one of several Bayesian tests that are reported. By
reducing the controversy regarding the specification of prior densities
on parameter values under individual hypotheses, UMPBTs can also be
used to focus attention on the specification of prior probabilities on
the hypotheses themselves. In the clinical trial example described
above, for example, the value of the success probability specified
under the alternative hypothesis may be less important in modeling
posterior model probabilities than incorporating information regarding
the outcomes of previous trials on related supplements. Such would be
the case if numerous previous trials of similar agents had failed to
provide evidence of increased treatment efficacy.

UMPBTs possess certain favorable properties not shared by other
objective Bayesian methods. For instance, most objective Bayesian tests
implicitly define local alternative prior densities on model parameters
under the alternative hypothesis [e.g., \citet{jeffreys39,ohagan95,berger96}]. As demonstrated in \citet{johnson10}, however, the use of
local alternative priors makes it difficult to accumulate evidence in
favor of a true null hypothesis. This means that many objective
Bayesian methods are only marginally better than classical significance
tests in summarizing evidence in favor of the null hypothesis. For
small to moderate sample sizes, UMPBTs produce alternative hypotheses
that correspond to nonlocal alternative prior densities, which means
that they are able to provide more balanced summaries of evidence
collected in favor of true null and true alternative hypotheses.

UMPBTs also possess certain unfavorable properties. Like many objective
Bayesian methods, UMPBTs can violate the likelihood principle, and
their behavior in large sample settings can lead to inconsistency if
evidence thresholds are held constant. And the alternative hypotheses
generated by UMPBTs are neither vague nor noninformative. Further
comments and discussion regarding these issues are provided below.

In order to define UMPBTs, it useful to first review basic properties
of Bayesian hypothesis tests.
In contrast to classical statistical hypothesis tests, Bayesian
hypothesis tests are based on comparisons of the posterior
probabilities assigned to competing hypotheses. In parametric tests,
competing hypotheses are characterized by the prior densities that they
impose on the parameters that define a sampling density shared by both
hypotheses. Such tests comprise the focus of this article.
Specifically, it is assumed throughout that the posterior odds between
two hypotheses $H_1$ and $H_0$ can be expressed as
%
%
\begin{equation}
\label{posteriorOdds} \frac{\bP(H_1 \con{\mathbf x})}{\bP(H_0 \con{\mathbf x})} = \frac
{m_1({\mathbf x})}{m_0({\mathbf x})} \times
\frac{p(H_1)}{p(H_0)},
\end{equation}
where
$
\mathrm{BF}_{10}({\mathbf x}) = m_1({\mathbf x})/m_0({\mathbf x})
$
is the Bayes factor between hypotheses $H_1$ and~$H_0$,
%
%
\begin{equation}
m_i({\mathbf x}) = \int_\Theta f({\mathbf x}\con\bft)
\pi_i(\bft\con H_i) \,d\bft
\end{equation}
is the marginal density of the data under hypothesis $H_i$, $f({\mathbf
x}\con\bft)$ is the sampling density of data ${\mathbf x}$ given $\bft$,
$\pi_i(\bft\con H_i)$ is the prior density on $\bft$ under $H_i$ and
$p(H_i)$ is the prior probability assigned to hypothesis $H_i$, for
$i=0,1$. The marginal prior density for $\bft$ is thus
\[
\pi(\theta) = \pi_0(\theta\con H_0)
P(H_0) + \pi_1(\theta\con H_1)
p(H_1).
\]
When there is no possibility of confusion, $\pi_i(\theta\con H_i)$
will be denoted more simply by $\pi_i(\theta)$. The parameter space
is denoted by $\Theta$ and the sample space by $\mathcal{X}$. The
logarithm of the Bayes factor is called the \emph{weight of evidence.}
All densities are assumed to be defined with respect to an appropriate
underlying measure (e.g., Lebesgue or counting measure).


Finally, assume that one hypothesis---the null hypothesis $H_0$---is
fixed on the basis of scientific considerations, and that the
difficulty in constructing a Bayesian hypothesis test arises from the
requirement to specify an alternative hypothesis. This assumption
mirrors the situation encountered in classical hypothesis tests in
which the null hypothesis is known, but no alternative hypothesis is
defined. In the clinical trial example, for instance, the null
hypothesis corresponds to the assumption that the success probability
of the new treatment equals that of the standard treatment, but there
is no obvious value (or prior probability density) that should be
assigned to the treatment's success probability under the alternative
hypothesis that it is better than the standard of care.

With these assumptions and definitions in place, it is worthwhile to
review a property of Bayes factors that pertains when the prior density
defining an alternative hypothesis is misspecified. Let $\pi_t(\bft
\con H_1)=\pi_t(\bft)$ denote the ``true'' prior density on $\bft$
under the assumption that the alternative hypothesis is true, and let
$m_t({\mathbf x})$ denote the resulting marginal density of the data. In
general $\pi_t(\bft)$ is not known, but it is still possible to
compare the properties of the weight of evidence that would be obtained
by using the true prior density under the alternative hypothesis to
those that would be obtained using some other prior density. From a
frequentist perspective, $\pi_t$ might represent a point mass
concentrated on the true, but unknown, data generating parameter. From
a Bayesian perspective, $\pi_t$ might represent a summary of existing
knowledge regarding $\bft$ before an experiment is conducted. Because
$\pi_t$ is not available, suppose that $\pi_1(\bft\con H_1)=\pi
_1(\bft)$ is instead used to represent the prior density, again under
the assumption that the alternative hypothesis is true. Then it follows
from Gibbs's inequality that
\begin{eqnarray*}
&&
\int_{\cal X} m_t({\mathbf x}) \log \biggl[
\frac{m_t({\mathbf x})}{m_0({\mathbf
x})} \biggr] \,d{\mathbf x} - \int_{\cal X}
m_t({\mathbf x}) \log \biggl[ \frac{m_1({\mathbf x})}{m_0({\mathbf
x})} \biggr] \,d{\mathbf x}
\\
&&\qquad= \int_{\cal X} m_t({\mathbf x}) \log \biggl[
\frac{m_t({\mathbf x})}{m_1({\mathbf
x})} \biggr] \,d{\mathbf x}
\\
&&\qquad\geq 0.
\end{eqnarray*}
That is,
%
%
\begin{equation}
\label{gibbs} \int_{\cal X} m_t({\mathbf x}) \log \biggl[
\frac{m_t({\mathbf x})}{m_0({\mathbf
x})} \biggr] \,d{\mathbf x} \geq \int_{\cal X}
m_t({\mathbf x}) \log \biggl[ \frac{m_1({\mathbf x})}{m_0({\mathbf
x})} \biggr] \,d{\mathbf x},
\end{equation}
which means that the expected weight of evidence in favor of the
alternative hypothesis is always decreased when $\pi_1(\bft)$ differs
from $\pi_t(\bft)$ (on a set with measure greater than 0). In
general, the UMPBTs described below will thus decrease the average
weight of evidence obtained in favor of a true alternative hypothesis.
\emph{In other words, the weight of evidence reported from a UMPBT will
tend to underestimate the actual weight of evidence provided by an
experiment in favor of a true alternative hypothesis.}

Like classical statistical hypothesis tests, the tangible consequence
of a Bayesian hypothesis test is often the rejection of one hypothesis,
say $H_0$, in favor of the second, say $H_1$. In a Bayesian test, the
null hypothesis is rejected if the posterior probability of $H_1$
exceeds a certain threshold. Given the prior odds between the
hypotheses, this is equivalent to determining a threshold, say $\gamma
$, over which the Bayes factor between $H_1$ and $H_0$ must fall in
order to reject $H_0$ in favor of $H_1$. It is therefore of some
practical interest to determine alternative hypotheses that maximize
the probability that the Bayes factor from a test exceeds a specified threshold.

With this motivation and notation in place, a UMPBT($\gamma$) may be
formally defined as follows.

\begin{definition*}
A uniformly most powerful Bayesian test for
evidence threshold $\gamma>0$ in favor of the alternative hypothesis
$H_1$ against a fixed null hypothesis~$H_0$, denoted by UMPBT($\gamma
$), is a Bayesian hypothesis test in which the Bayes factor for the
test satisfies the following inequality for any $\bft_t \in\Theta$
and for all alternative hypotheses $H_2\dvtx  \bft\sim\pi_2(\bft)$:
%
%
\begin{equation}
\bP_{\bft_t} \bigl[ \mathrm{BF}_{10}({\mathbf x}) > \gamma \bigr] \geq\bP
_{\bft_t} \bigl[ \mathrm{BF}_{20}({\mathbf x}) > \gamma \bigr].
\end{equation}
\end{definition*}

In other words, the UMPBT($\gamma$) is a Bayesian test for which the
alternative hypothesis is specified so as to maximize the probability
that the Bayes factor $\mathrm{BF}_{10}({\mathbf x})$ exceeds the evidence threshold
$\gamma$ for all possible values of the data generating parameter
$\bft_t $.

The remainder of this article is organized as follows. In the next
section, UMPBTs are described for one-parameter exponential family
models. As in the case of UMPTs, a general prescription for
constructing UMPBTs is available only within this class of densities.
Specific techniques for defining UMPBTs or approximate UMPBTs outside
of this class are described later in Sections~\ref{extensions} and~\ref{sec5}.
In applying UMPBTs to one parameter exponential family models, an
approximate equivalence between type I errors for UMPTs and the Bayes
factors obtained from UMPBTs is exposed.

In Section~\ref{sec3}, UMPBTs are applied in two canonical testing situations:
the test of a binomial proportion, and the test of a normal mean. These
two tests are perhaps the most common tests used by practitioners of
statistics. The binomial test is illustrated in the context of a
clinical trial, while the normal mean test is applied to evaluate
evidence reported in support of the Higgs boson. Section~\ref{extensions} describes several settings outside of one parameter
exponential family models for which UMPBTs exist. These include cases
in which the nuisance parameters under the null and alternative
hypothesis can be considered to be equal (though unknown), and
situations in which it is possible to marginalize over nuisance
parameters to obtain expressions for data densities that are similar to
those obtained in one-parameter exponential family models. Section~\ref{sec5}
describes approximations to UMPBTs obtained by specifying alternative
hypotheses that depend on data through statistics that are ancillary to
the parameter of interest. Concluding comments appear in Section~\ref{discussion}.

\section{One-parameter exponential family models}\label{sec2}

Assume that $\{x_1,\ldots,\break x_n\} \equiv{\mathbf x}$ are i.i.d. with a
sampling density (or probability mass function in the case of discrete
data) of the form
%
%
\begin{equation}
\label{ef1} f(x \con\theta) = h(x) \exp \bigl[ \eta(\theta) T(x) - A(\theta)
\bigr],
\end{equation}
where $T(x)$, $h(x)$, $\eta(\theta)$ and $A(\theta)$ are known
functions, and $\eta(\theta)$ is monotonic. Consider a one-sided test
of a point null hypothesis $H_0\dvtx  \theta=\theta_0$ against an
arbitrary alternative hypothesis. Let $\gamma$ denote the evidence
threshold for a UMPBT($\gamma)$, and assume that the value of $\theta
_0$ is fixed.

\begin{lemm}\label{lm1}
Assume the conditions of the previous paragraph pertain, and define
$g_\gamma(\theta,\theta_0)$ according to
%
%
\begin{equation}
\label{gdef} g_\gamma(\theta,\theta_0) = \frac{\log(\gamma) + n[ A(\theta
)-A(\theta_0) ]}{\eta(\theta)-\eta(\theta_0)}.
\end{equation}
In addition, define $u$ to be 1 or $-1$ according to whether $\eta
(\theta)$ is monotonically increasing or decreasing, respectively, and
define $v$ to be either 1 or $-1$ according to whether the alternative
hypothesis requires $\theta$ to be greater than or less than~$\theta
_0$, respectively. Then a $\operatorname{UMPBT}(\gamma)$ can be obtained by
restricting the support of $\pi_1(\theta)$ to values of $\theta$
that belong to the set
%
%
\begin{equation}
\label{obj} \argmin_{\theta} uv g_\gamma(\theta,
\theta_0).
\end{equation}
\end{lemm}

\begin{pf}
Consider the case in which the alternative
hypothesis requires $\theta$ to be greater than $\theta_0$ and $\eta
(\theta)$ is increasing (so that $uv=1$), and let $\theta_t$ denote
the true (i.e., data-generating) parameter for ${\mathbf x}$ under (\ref
{ef1}). Consider first simple alternatives for which the prior on
$\theta$ is a point mass at $\theta_1$. Then
%
%
\begin{eqnarray}\label{ineq}
\bP_{\theta_t} ( \mathrm{BF}_{10} > \gamma ) &=& \bP_{\theta
_t} \bigl[
\log(\mathrm{BF}_{10}) > \log(\gamma) \bigr]
\nonumber\\[-8pt]\\[-8pt]
&=& \bP_{\theta_t} \Biggl\{ \sum_{i=1}^n
T(x_i) > \frac{\log(\gamma
) + n[A(\theta_1)-A(\theta_0)]}{\eta(\theta_1)-\eta(\theta_0)}\Biggr\}.
\nonumber
\end{eqnarray}
It follows that the probability in (\ref{ineq}) achieves its maximum
value when the right-hand side of the inequality is minimized,
regardless of the distribution of $\sum T(x_i)$.\vadjust{\goodbreak}

Now consider composite alternative hypotheses, and define an indicator
function $s$ according to
%
%
\begin{equation}
s({\mathbf x},\theta) = \operatorname{Ind} \Biggl( \exp \Biggl\{ \bigl[\eta(\theta)-\eta (
\theta_0)\bigr] \sum_{i=1}^n
T(x_i) -n\bigl[A(\theta)-A(\theta_0)\bigr] \Biggr\} >
\gamma \Biggr).\hspace*{-28pt}
\end{equation}
Let $\theta^*$ be a value that minimizes $g_\gamma(\theta,\theta
_0)$. Then it follows from (\ref{ineq}) that
%
%
\begin{equation}
\label{dvec1}
s({\mathbf x},\theta)\leq s\bigl({\mathbf x},\theta^*\bigr) \qquad\mbox{for all }
{\mathbf x}.
\end{equation}
This implies that
%
%
\begin{equation}
\int_\Theta s({\mathbf x},\theta) \pi(\theta) \,d\theta\leq s\bigl({
\mathbf x},\theta^*\bigr)
\end{equation}
for all probability densities $\pi(\theta)$. It follows that
%
%
\begin{equation}
\label{dvec2} \bP_{\theta_t} ( \mathrm{BF}_{10} > \gamma ) = \int
_{\mathcal X} s({\mathbf x},\theta) f({\mathbf x}\con\theta_t) \,dx
\end{equation}
is maximized by a prior that concentrates its mass on the set for which
$g_\gamma(\theta,\theta_0)$ is minimized.

The proof for other values of $(u,v)$ follows by noting that the
direction of the inequality in (\ref{ineq}) changes according to the
sign of $\eta(\theta_1)-\eta(\theta_0)$.
\end{pf}

It should be noted that in some cases the values of $\theta$ that
maximize $\bP_{\theta_t}  ( \mathrm{BF}_{10} > \gamma )$ are not
unique. This might happen if, for instance, no value of the sufficient
statistic obtained from the experiment could produce a Bayes factor
that exceeded the $\gamma$ threshold. For example, it would not be
possible to obtain a Bayes factor of 10 against a null hypothesis that
a binomial success probability was 0.5 based on a sample of size $n=1$.
In that case, the probability of exceeding the threshold is 0 for all
values of the success probability, and a unique UMPBT does not exist.
More generally, if $T(x)$ is discrete, then many values of $\theta_1$
might produce equivalent tests. An illustration of this phenomenon is
provided in the first example.

\subsection{Large sample properties of UMPBTs}\label{sec21}\label{largesample}
Asymptotic properties of\break UMPBTs can most easily be examined for tests
of point null hypotheses for a canonical parameter in one-parameter
exponential families. Two properties of UMPBTs in this setting are
described in the following lemma.
%
\begin{lemm}\label{lem2}
Let $X_1,\ldots,X_n$ represent a random sample drawn from a density
expressible in the form (\ref{ef1}) with $\eta(\theta)=\theta$, and
consider a test of the precise null hypothesis $H_0\dvtx  \theta=\theta
_0$. Suppose that $A(\theta)$ has three bounded derivatives in a
neighborhood of $\theta_0$, and let $\theta^*$ denote a value of
$\theta$ that defines a $\operatorname{UMPBT}(\gamma)$ test and satisfies
%
%
\begin{equation}
\label{dg} \frac{dg_\gamma(\theta^*,\theta_0)}{d\theta} = 0.\vadjust{\goodbreak}
\end{equation}
Then the following statements apply:
\begin{longlist}[(2)]
\item[(1)] For some $t \in(\theta_0,\theta^*)$,
%
%
\begin{equation}
\label{pitman} \bigl\llvert \theta^*-\theta_0 \bigr\rrvert = \sqrt{
\frac{2\log(\gamma
)}{nA''(t)} }.
\end{equation}
\item[(2)] Under the null hypothesis,
%
%
\begin{equation}
\label{clt} \log(\mathrm{BF}_{10}) \rightarrow N \bigl( -\log(\gamma), 2\log(
\gamma) \bigr) \qquad\mbox{as } n \rightarrow\infty.
\end{equation}
\end{longlist}
\end{lemm}

\begin{pf}
The first statement follows immediately from
(\ref{dg}) by expanding $A(\theta)$ in a Taylor series around $\theta^{*}$.
The second statement follows by noting that the weight of evidence can
be expressed as
\[
\log(\mathrm{BF}_{10}) = \bigl(\theta^*-\theta_0\bigr) \sum
_{i=1}^n T(x_i) - n\bigl[A
\bigl(\theta^*\bigr)-A(\theta_0)\bigr].
\]
Expanding in a Taylor series around $\theta_0$ leads to
%
%
\begin{equation}
\label{bigN} \log(\mathrm{BF}_{10}) = \sqrt{\frac{2\log(\gamma)}{nA''(t)}} \Biggl[ \sum
_{i=1}^n T(x_i) -
nA'(\theta_0) - \frac{n}{2} A''(
\theta_0) \sqrt {\frac{2\log(\gamma)}{nA''(t)}} \Biggr] +
\varepsilon,\hspace*{-28pt}
\end{equation}
where $\varepsilon$ represents a term of order $O(n^{-1/2})$. From
properties of exponential family models, it is known that
\[
E_{\theta_0} \bigl[T(x_i)\bigr] = A'(
\theta_0) \quad\mbox{and}\quad \operatorname{Var}_{\theta_0}\bigl(T(x_i)
\bigr) = A''(\theta_0).
\]
Because $A(\theta)$ has three bounded derivatives in a neighborhood of
$\theta_0$, $[A'(t)-A'(\theta_0)]$ and $[A''(t)-A''(\theta_0)]$ are
order $O(n^{-1/2})$, and the statement follows by application of the
central limit theorem.
\end{pf}

Equation (\ref{pitman}) shows that the difference $|\theta^*-\theta
_0|$ is $O(n^{-1/2})$ when the evidence threshold $\gamma$ is held
constant as a function of $n$. In classical terms, this implies that
alternative hypotheses defined by UMPBTs represent Pitman sequences of
local alternatives [\citet{pitman49}]. This fact, in conjunction with
(\ref{clt}), exposes an interesting behavior of UMPBTs in large sample
settings, particularly when viewed from the context of the
Jeffreys--Lindley paradox [\citet{jeffreys39,lindley57}; see also \citet
{robert09}].

The Jeffreys--Lindley paradox (JLP) arises as an incongruity between
Bayesian and classical hypothesis tests of a point null hypothesis. To
understand the paradox, suppose that the prior distribution for a
parameter of interest under the alternative hypothesis is uniform on an
interval $I$ containing the null value $\theta_0$, and that the prior
probability assigned to the null hypothesis is $\pi_1$. If $\pi_1$ is
bounded away from~0, then it is possible for the null hypothesis to be
rejected in an $\alpha$-level significance test even when the
posterior probability assigned to the null hypothesis exceeds $1-\alpha
$. Thus, the anomalous behavior exhibited in the example of Table~\ref{tabl1}, in
which the null hypothesis was rejected in a significance test while
being supported by the data, is characteristic of a more general
phenomenon that may occur even in large sample settings.
To see that the null hypothesis can be rejected even when the posterior
odds are in its favor, note that for sufficiently large $n$ the width
of $I$ will be large relative to the posterior standard deviation of
$\theta$ under the alternative hypothesis. Data that are not ``too
far'' from $f_{\theta_0}$ may therefore be much more likely to have
arisen from the null hypothesis than from a density $f_\theta$ when
$\theta$ is drawn uniformly from $I$. At the same time, the value of
the test statistic based on the data may appear extreme given that
$f_{\theta_0}$ pertains.

For moderate values of $\gamma$, the second statement in Lemma~\ref{lem2} shows
that the weight of evidence obtained from a UMPBT is unlikely to
provide strong evidence in favor of either hypothesis when the null
hypothesis is true. When $\gamma=4$, for instance, an approximate 95\%
confidence interval for the weight of evidence extends only between
$(-4.65,1.88)$, no matter how large $n$ is. Thus, the posterior
probability of the null hypothesis does not converge to 1 as the sample
size grows. The null hypothesis is never fully accepted---nor the
alternative rejected---when the evidence threshold is held constant as
$n$ increases.

This large sample behavior of UMPBTs with fixed evidence thresholds is,
in a certain sense, similar to the JLP. When the null hypothesis is
true and $n$ is large, the probability of rejecting the null hypothesis
at a fixed level of significance remains constant at the specified
level of significance. For instance, the null hypothesis is rejected
5\% of the time in a standard 5\% significance test when the null
hypothesis is true, regardless of how large the sample size is.
Similarly, when $\gamma=4$, the probability that the weight of
evidence in favor of the alternative hypothesis will be greater than 0
converges to 0.20 as $n$ becomes large. Like the significance test,
there remains a nonzero probability that the alternative hypothesis
will be favored by the UMPBT even when the null hypothesis is true,
regardless of how large $n$ is.

On a related note, \citet{rousseau06} has demonstrated that a point
null hypothesis may be used as a convenient mathematical approximation
to interval hypotheses of the form $(\theta_0-\varepsilon,\theta
_0+\varepsilon)$ if $\varepsilon$ is sufficiently small. Her results
suggest that such an approximation is valid only if $\varepsilon< o(n)$.
The fact that UMPBT alternatives decrease at a rate of $O(n^{-1/2})$
suggests that UMPBTs may be used to test small interval hypotheses
around $\theta_0$, provided that the width of the interval satisfies
the constraints provided by Rousseau.

Further comments regarding the asymptotic properties of UMPBTs appear
in the discussion section.

\section{Examples}\label{sec3}\label{examples}

Tests of simple hypotheses in one-parameter exponential family models
continue to be the most common statistical hypothesis tests used by
practitioners. These tests play a central role in many science,
technology and business applications. In addition, the distributions of
many test statistics are asymptotically distributed as standard normal
deviates, which means that UMPBTs can be applied to obtain Bayes
factors based on test statistics [\citet{johnson05}]. This section
illustrates the use of UMPBT tests in two archetypical examples; the
first involves the test of a binomial success probability, and the
second the test of the value of a parameter estimate that is assumed to
be approximately normally distributed.

\subsection{Test of binomial success probability}\label{sec31}

Suppose $x \sim \operatorname{Bin}(n,p)$, and consider the test of a null hypothesis
$H_0\dvtx  p=p_0$ versus an alternative hypothesis $H_1\dvtx  p>p_0$. Assume that
an evidence threshold of $\gamma$ is desired for the test; that is,
the alternative hypothesis is accepted if $\mathrm{BF}_{10}>\gamma$.

From Lemma~\ref{lm1}, the UMPBT($\gamma$) is defined by finding $p_1$ that
satisfies $p_1>p_0$ and
%
%
\begin{equation}
\label{prop1} p_1 = \argmin_{p} \frac{\log(\gamma) - n[ \log(1-p)-\log
(1-p_0)]}{\log[p/(1-p)] - \log[ p_0/(1-p_0)]}.
\end{equation}
Although this equation cannot be solved in closed form, its solution
can be found easily using optimization functions available in most
statistical programs.

\subsubsection{Phase \textup{II} clinical trials with binary outcomes}\label{sec311}

To illustrate the resulting test in a real-world application that
involves small sample sizes, consider a one-arm Phase II trial of a new
drug intended to improve the response rate to a disease from the
standard-of-care rate of $p_0=0.3$. Suppose also that budget and time
constraints limit the number of patients that can be accrued in the
trial to $n=10$, and suppose that the new drug will be pursued only if
the odds that it offers an improvement over the standard of care are at
least 3:1. Taking $\gamma=3$, it follows from (\ref{prop1}) that
the UMPBT alternative is defined by taking $H_1\dvtx  p_1=0.525$. At this
value of $p_1$, the Bayes factor $\mathrm{BF}_{10}$ in favor of $H_1$ exceeds 3
whenever $6$ or more of the 10 patients enrolled in the trial respond
to treatment.

\begin{figure}

\includegraphics{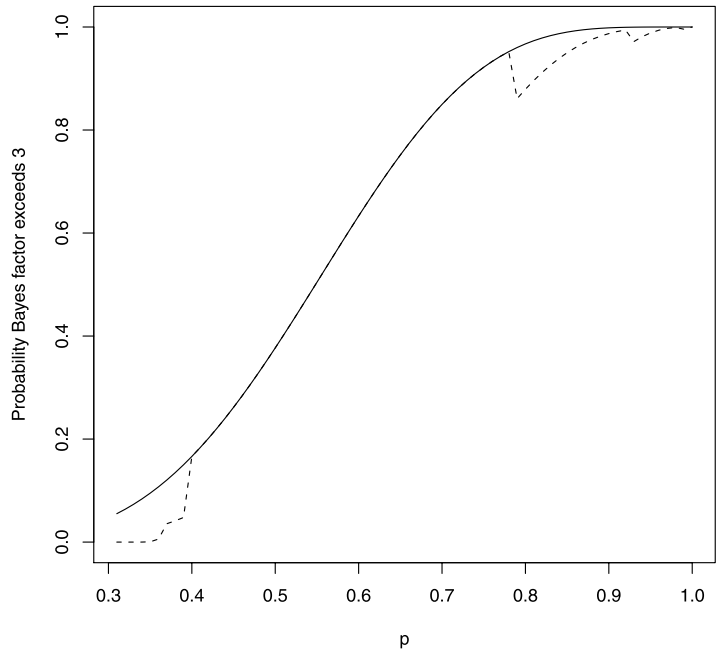}

\caption{Probability that the Bayes factor exceeds 3 plotted against
the data-generating parameter. The solid curve shows the probability of
exceeding 3 for the UMPBT. The dashed curve displays this probability
when the Bayes factor is calculated using the data-generating parameter.}
\label{binplot}
\end{figure}

A plot of the probability that $\mathrm{BF}_{10}$ exceeds $3$ as function of the
true response rate $p$ appears in Figure~\ref{binplot}. For
comparison, also plotted in this figure (dashed curve) is the
probability that $\mathrm{BF}_{10}$ exceeds 3 when $p_1$ is set to the
data-generating parameter, that is, when $p_1=p_t$.

Figure~\ref{binplot} shows that the probability that $\mathrm{BF}_{10}$ exceeds
3 when calculated under the true alternative hypothesis is
significantly smaller than it is under the UMPBT alternative for values
of $p<0.4$ and for values of $p>0.78$. Indeed, for values of $p<0.334$,
there is no chance that $\mathrm{BF}_{10}$ will exceed~3. This is so because
$(0.334/0.30)^x$ remains less than 3.0 for all $x\leq10$. The decrease
in the probability that the Bayes factor exceeds 3 for large values of
$p$ stems from the relatively small probability that these models
assign to the observation of intermediate values of $x$. For example,
when $p=0.8$, the probability of observing 6 out 10 successes is only
0.088, while the corresponding probability under $H_0$ is 0.037. Thus
$\mathrm{BF}_{10}=2.39$, and the evidence in favor of the true success
probability does not exceed 3. That is, the discontinuity in the dashed
curve at $p\approx0.7$ occurs because the Bayes factor for this test
is not greater than 3 when $x=6$. Similarly, the other discontinuities
in the dashed curve occur when the rejection region for the Bayesian
test (i.e., values of $x$ for which the Bayes factor is greater than 3)
excludes another immediate value of $x$. The dashed and solid curves
agree for all Bayesian tests that produce Bayes factors that exceed 3
for all values of $x\geq6$.

It is also interesting to note that the solid curve depicted in
Figure~\ref{binplot} represents the power curve for an approximate 5\%
one-sided significance test of the null hypothesis that $p=0.3$ [note
that $P_{0.3}(X\geq6) = 0.047$]. This rejection region for the 5\%
significance test also corresponds to the region for which the Bayes
factor corresponding to the UMPBT($\gamma$) exceeds $\gamma$ for all
values of $\gamma\in(2.36,6.82)$. If equal prior probabilities are
assigned to $H_0$ and $H_1$, this suggests that a $p$-value of 0.05 for
this test corresponds roughly to the assignment of a posterior
probability between $(1.0/7.82,1.0/3.36) = (0.13,0.30)$ to the null
hypothesis. This range of values for the posterior probability of the
null hypothesis is in approximate agreement with values suggested by
other authors, for example, \citet{berger87}.

This example also indicates that a UMPBT can result in large type I
errors if the threshold $\gamma$ is chosen to be too small. For
instance, taking $\gamma=2$ in this example would lead to type I
errors that were larger than 0.05.

\begin{figure}

\includegraphics{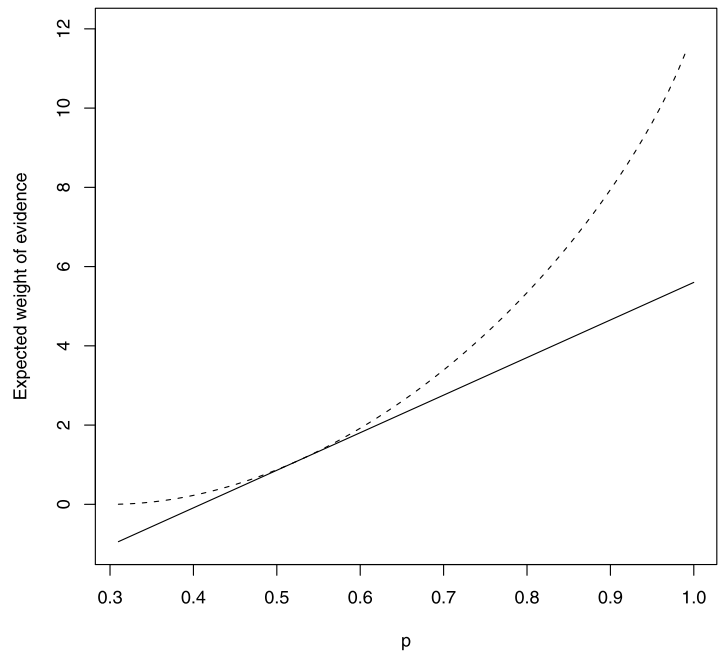}

\caption{Expected weight of evidence produced by a
$\operatorname{UMPBT}(\gamma)$
against a null hypothesis that $p_0=0.3$ when the sample size is $n=10$
(solid curve), versus the expected weight of evidence observed using
the data-generating success probability at the alternative hypothesis
(dashed curve). The data-generating parameter value is displayed on the
horizontal axis.}
\label{bininfo}
\end{figure}

It is important to note that the UMPBT does not provide a test that
maximizes the expected weight of evidence, as equation (\ref{gibbs})
demonstrates. This point is illustrated in Figure~\ref{bininfo}, which
depicts the expected weight of evidence obtained in favor of $H_1$ by a
solid curve as the data-generating success probability is varied in
%
$(0.3,1.0)$. For comparison, the dashed curve shows the expected weight
of evidence obtained as a function of the true parameter value. As
predicted by the inequality in (\ref{gibbs}), on average the UMPBT
provides less evidence in favor of the true alternative hypothesis for
all values of $p\in(0.3,1.0)$ except $p=0.525$, the UMPBT value.

\subsection{\texorpdfstring{Test of normal mean, $\sigma^2$ known}
{Test of normal mean, sigma2 known}}\label{sec32}
Suppose $x_i$, $i=1,\ldots,n$ are i.i.d. $N(\mu,\sigma^2)$ with
$\sigma^2$ known. The null hypothesis is $H_0\dvtx  \mu=\mu_0$, and the
alternative hypothesis is accepted if $\mathrm{BF}_{10}>\gamma$. Assuming that
the alternative hypothesis takes the form $H_1\dvtx  \mu=\mu_1$ in a
one-sided test, it follows that
%
%
\begin{equation}
\label{normalBF} \log(\mathrm{BF}_{10}) = \frac{n}{\sigma^2} \biggl[ \bar{{\mathbf
x}}(\mu_1-\mu _0) + \frac{1}{2}\bigl(
\mu_0^2-\mu_1^2\bigr) \biggr].
\end{equation}
If the data-generating parameter is $\mu_t$, the probability that
$\mathrm{BF}_{10}$ is greater than $\gamma$ can be written as
%
%
\begin{equation}
\label{normalbf} \bP_{\mu_t} \biggl[ (\mu_1-
\mu_0) \bar{{\mathbf x}} > \frac{\sigma^2
\log(\gamma)}{n} - \frac{1}{2} \bigl(
\mu_0^2-\mu_1^2\bigr) \biggr].
\end{equation}
If $\mu_1>\mu_0$, then the UMPBT($\gamma$) value of $\mu_1$ satisfies
%
%
\begin{equation}
\argmin_{\mu_1} \frac{\sigma^2 \log(\gamma)}{n(\mu_1-\mu_0)} +\frac{1}{2}(
\mu_0+\mu_1).
\end{equation}
Conversely, if $\mu_1<\mu_0$, then optimal value of $\mu_1$ satisfies
%
%
\begin{equation}
\argmin_{\mu_1} \frac{\sigma^2 \log(\gamma)}{n(\mu_1-\mu_0)} +\frac{1}{2}(\mu
_0+\mu_1).
\end{equation}
It follows that the UMPBT($\gamma$) value for $\mu_1$ is given by
%
%
\begin{equation}
\label{normal} \mu_1 = \mu_0 \pm\sigma\sqrt{
\frac{2 \log\gamma}{n}},
\end{equation}
depending on whether $\mu_1>\mu_0$ or $\mu_1<\mu_0$.

\begin{figure}

\includegraphics{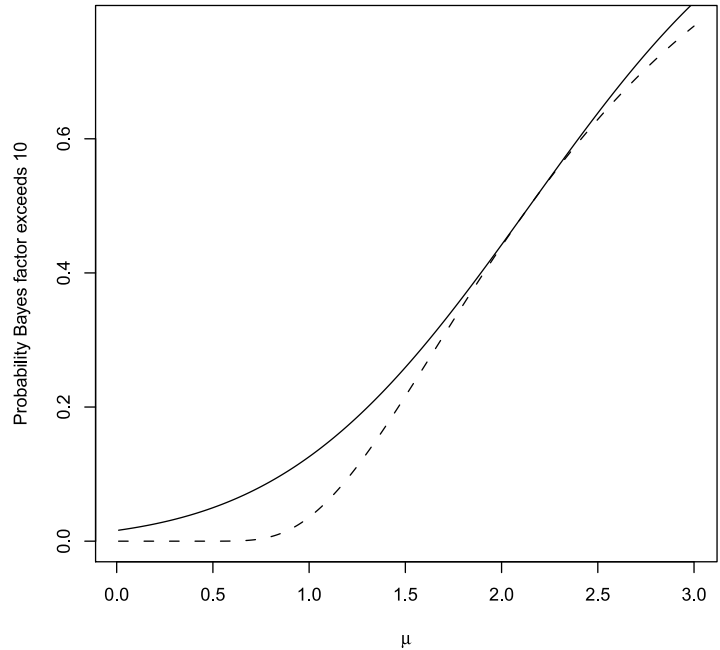}

\caption{Probability that Bayes factor in favor of UMPBT alternative
exceeds 10 when $\mu_0=0$ and $n=1$ (solid curve). The dashed curve
displays this probability when the Bayes factor is calculated under the
alternative hypothesis that $\mu_1$ equals the data-generating
parameter (displayed on the horizontal axis).}
\label{muplot}
\end{figure}

Figure~\ref{muplot} depicts the probability that the Bayes factor
exceeds $\gamma=10$ when testing a null hypothesis that $\mu=0$ based
on a single, standard normal observation (i.e., $n=1$, $\sigma^2=1$).
In this case, the UMPBT(10) is obtained by taking $\mu_1=2.146$. For
comparison, the probability that the Bayes factor exceeds 10 when the
alternative is defined to be the data-generating parameter is depicted
by the dashed curve in the plot.

UMPBTs can also be used to interpret the evidence obtained from
classical UMPTs. In a classical one-sided test of a normal mean with
known variance, the null hypothesis is rejected if
\[
\bar{x} > \mu_0 + z_{\alpha} \frac{\sigma}{\sqrt{n}},
\]
where $\alpha$ is the level of the test designed to detect $\mu_1>\mu
_0$. In the UMPBT, from (\ref{normalBF})--(\ref{normalbf}) it follows
that the null hypothesis is rejected if
\[
\bar{x} > \frac{\sigma^2 \log(\gamma)}{n(\mu_1-\mu_0)} + \frac
{1}{2}(\mu_1+
\mu_0).
\]
Setting $\mu_1=\mu_0+\sigma\sqrt{2\log(\gamma)/n}$ and equating
the two rejection regions, it follows that the rejection regions for
the two tests are identical if
%
%
\begin{equation}
\label{gammaequiv} \gamma= \exp \biggl( \frac{z_{\alpha}^2}{2} \biggr).
\end{equation}
For the case of normally distributed data, it follows that
%
%
\begin{equation}
\label{rr} \mu_1 = \mu_0 + \frac{\sigma}{\sqrt{n}}
z_{\alpha},
\end{equation}
which means that the alternative hypothesis places $\mu_1$ at the
boundary of the UMPT rejection region.

The close connection between the UMPBT and UMPT for a normal mean makes
it relatively straightforward to examine the relationship between the
$p$-values reported from a classical test and either the Bayes factors
or posterior probabilities obtained from a Bayesian test. For example,
significance tests for normal means are often conducted at the 5\%
level. Given this threshold of evidence for rejection of the null
hypothesis, the one-sided $\gamma$ threshold corresponding to the 5\%
significance test is 3.87, and the UMPBT alternative is $\mu_1= \mu
_0+1.645 \sigma/\sqrt{n}$. If we assume that equal prior
%
\begin{figure}

\includegraphics{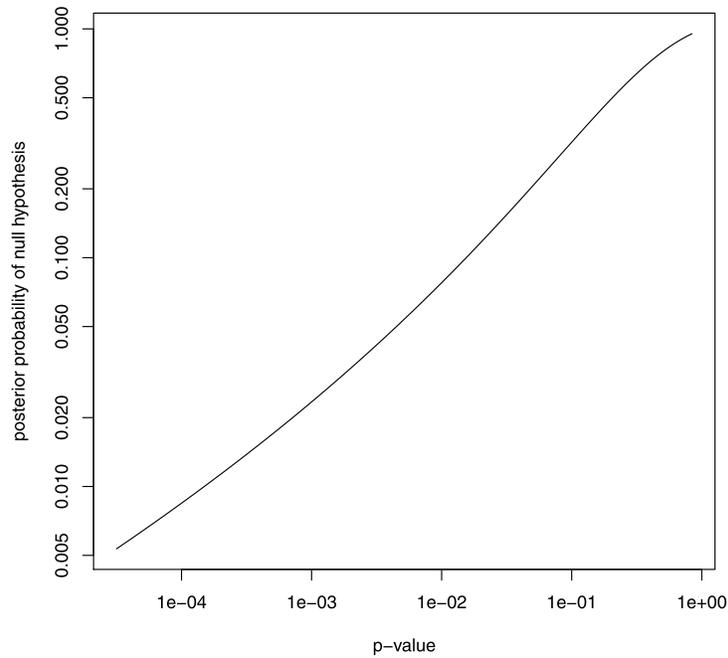}

\caption{Correspondence between $p$-values and posterior model
probabilities for a UMPBT test derived from a 5\% test. This plot
assumes equal prior probabilities were assigned to the null and
alternative hypotheses. Note that both axes are displayed on the
logarithmic scale.}
\label{correspondence}
\end{figure}
probabilities are assigned to the null and alternative hypotheses, then
a correspondence between $p$-values and posterior probabilities
assigned to the null hypothesis is easy to establish. This
correspondence is depicted in Figure~\ref{correspondence}. For
instance, this figure shows that a $p$-value of 0.01 corresponds to the
assignment of posterior probability 0.08 to the null hypothesis.

\subsubsection{Evaluating evidence for the Higgs boson}\label{sec321}
On July 4, 2012, scientists at CERN made the following announcement:

\begin{quote}
We observe in our data clear signs of a new particle, at the level of 5
sigma, in the mass region around 126 gigaelectronvolts (GeV).
(\href{http://press.web.cern.ch/press/PressReleases/Releases2012/PR17.12E.html}{http://press.web.cern.ch/}
\href{http://press.web.cern.ch/press/PressReleases/Releases2012/PR17.12E.html}{press/PressReleases/Releases2012/PR17.12E.html}).
\end{quote}

\noindent In very simplified terms, the 5 sigma claim can be explained by fitting
a model for a Poisson mean that had the following approximate form:
\[
\mu(x) = \exp\bigl(a_0 + a_1 x + a_2
x^2\bigr) + s \phi(x;m, w).
\]
Here,
$x$ denotes mass in GeV, $\{a_i\}$ denote nuisance parameters that
model background events, $s$ denotes signal above background, $m$
denotes the mass of a new particle, $w$ denotes a convolution parameter
and $\phi(x;m,w)$ denotes a Gaussian density centered on $m$ with
standard deviation $w$ [\citet{prosper12}]. Poisson events
collected from a series of high energy experiments conducted in the
Large Hadron Collider (LHC) at CERN provide the data to estimate the
parameters in this stylized model. The background parameters $\{ a_i\}$
are considered nuisance parameters. Interest, of course, focuses on
testing whether $s>0$ at a mass location $m$. The null hypothesis is
that $s=0$ for all $m$.

The accepted criterion for declaring the discovery of a new particle in
the field of particle physics is the 5 sigma rule, which in this case
requires that the estimate of $s$ be 5 standard errors from 0
(\href{http://public.web.cern.ch/public/}{http://public.web.cern.ch/public/}).

Calculation of a Bayes factor based on the original mass spectrum data
is complicated by the fact that prior distributions for the nuisance
parameters $\{a_i\}$, $m$, and $w$ are either not available or are not
agreed upon. For this reason, it is more straightforward to compute a
Bayes factor for these data based on the test statistic $z=\hat
{s}/\operatorname{se}(\hat{s})$ where $\hat{s}$ denotes the maximum likelihood
estimate of $s$ and $\operatorname{se}(\hat{s})$ its standard error [Johnson
(\citeyear{johnson05,johnson08})].
To perform this test, assume that under the null hypothesis $z$ has a
standard normal distribution, and that under the alternative hypothesis
$z$ has a normal distribution with mean $\mu$ and variance 1.

In this context, the 5 sigma rule for declaring a new particle
discovery means that a new discovery can only be declared if the test
statistic $z>5$. Using equation (\ref{gammaequiv}) to match the
rejection region of the classical significance test to a UMPBT($\gamma
$) implies that the corresponding evidence threshold is $\gamma=
\exp(12.5)\approx27\mbox{,}000$. In other words, a
Bayes factor of approximately $\gamma=\exp(12.5)\approx27\mbox{,}000$
corresponds to the 5 sigma rule required to accept the alternative
hypothesis that a new particle has been found.

It follows from the discussion following equation (\ref{rr}) that the
alternative hypothesis for the UMPBT alternative is $\mu_1=5$. This
value is calculated under the assumption that the test statistic $z$
has a standard normal distribution under the null hypothesis [i.e.,
$\sigma=1$ and $n=1$ in (\ref{normal})]. If the observed value of $z$
was exactly 5, then the Bayes factor in favor of a new particle would
be approximately 27,000. If the observed value was, say 5.1, then the
Bayes factor would be $\exp(-0.5[0.1^2-5.1^2])=44\mbox{,}000$. These values
suggest very strong evidence in favor of a new particle, but perhaps
not as much evidence as might be inferred by nonstatisticians by the
report of a $p$-value of $3\times10^{-7}$.

There are, of course, a number of important caveats that should be
considered when interpreting the outcome of this analysis. This
analysis assumes that an experiment with a fixed endpoint was
conducted, and that the UMPBT value of the Poisson rate at 126 GeV was
of physical significance. Referring to (\ref{normal}) and noting that
the asymptotic standard error of $z$ decreases at rate $\sqrt{n}$, it
follows that the UMPBT alternative hypothesis favored by this analysis\vadjust{\goodbreak}
is $O(n^{-1/2})$. For sufficiently large $n$, systematic errors in the
estimation of the background rate could eventually lead to the
rejection of the null hypothesis in favor of the hypothesis of a new
particle. This is of particular concern if the high energy experiments
were continued until a 5 sigma result was obtained. Further comments
regarding this point appear in the discussion section.

\subsection{Other one-parameter exponential family models}\label{sec33}
Table~\ref{objective} provides the functions that must be minimized to
obtain UMPBTs for a number of common exponential family models. The
objective functions listed in this table correspond to the function
%
\begin{table}
\caption{Common one parameter exponential family models for which
$\operatorname{UMPBT}(\gamma)$ exist}
\label{objective}
\begin{tabular*}{\tablewidth}{@{\extracolsep{\fill}}lcc@{}}
\hline
\textbf{Model} & \textbf{Test} & \textbf{Objective function} \\
\hline
Binomial & $p_1>p_0$ & $ \{\log(\gamma) - n \log
[(1-p)/(1-p_0)] \} (\log\{[p(1-p_0)]/[(1-p)p_0]\}
)^{-1}$ \\[2pt]
Exponential & $\mu_1 >\mu_0 $ &
$  \{ \log(\gamma) + n [\log(\mu_1)-\log(\mu_0)] \}
[ 1/\mu_0 - 1/\mu_1 ]^{-1}$ \\[2pt]
Neg. Bin.& $p_1>p_0$ &
$  \{\log(\gamma) - r\log[(1-p_1)/(1-p_0)] \} [\log
(p_1) - \log(p_0) ]^{-1} $
\\[2pt]
Normal & $\sigma_1^2 >\sigma_0^2 $ & $
\{ 2\sigma_1^2\sigma_0^2  ( \log(\gamma)+\frac{n}{2}
[\log(\sigma_1^2)-\log(\sigma_0^2)]  ) \} [ \sigma
_1^2-\sigma_0^2 ]^{-1}$ \\[2pt]
Normal& $\mu_1 >\mu_0 $ & $  [ \sigma^2 \log(\gamma)
] ( \mu_1-\mu_0 )^{-1} + \frac{1}{2} (\mu_0+\mu_1) $ \\[2pt]
Poisson & $\mu_1 >\mu_0 $ & $ [ \log(\gamma) + n(\mu_1-\mu
_0) ]  [ \log(\mu_1) - \log(\mu_0) ]^{-1}$ \\
\hline
\end{tabular*}
\end{table}
$g_\gamma(\cdot,\cdot)$ specified in Lemma~\ref{lm1} with $v=1$. The
negative binomial is parameterized by the fixed number of failures $r$
and random number of successes $x=0,1,\ldots$ observed before
termination of sampling. The other models are parameterized so that
$\mu$ and $p$ denote means and proportions, respectively, while
$\sigma^2$ values refer to variances.

\section{Extensions to other models}\label{sec4}\label{extensions}
Like UMPTs, UMPBTs are most easily defined within one-parameter
exponential family models. In unusual cases, UMPBTs can be defined for
data modeled by vector-valued exponential family models, but in general
such extensions appear to require stringent constraints on nuisance parameters.

One special case in which UMPBTs can be defined for a $d$-dimensional
parameter $\bft$ occurs when the density of an observation can be
expressed as
%
%
\begin{equation}
\label{vector} f(x\con\bft) = h(x) \exp \Biggl[ \sum
_{i=1}^d \eta_i(\bft)
T_i(x) - A(\bft) \Biggr],
\end{equation}
and all but one of the $\eta_i(\bft)$ are constrained to have
identical values under both hypotheses. To understand how a UMPBT can
be defined in this case, without loss of generality suppose that $\eta
_i(\bft)$, $i=2,\ldots,d$ are constrained to have the same value under
both the null and alternative hypotheses, and that the null and
alternative hypotheses are defined by $H_0\dvtx  \theta_1 = \theta_{1,0}$
and $H_1\dvtx  \theta_1>\theta_{1,0}$. For simplicity, suppose further
that $\eta_1$ is a monotonically increasing function.\vadjust{\goodbreak}

As in Lemma~\ref{lm1}, consider first simple alternative hypotheses expressible
as $H_1\dvtx  \theta_1 = \theta_{1,1}$. Let $\bft_0 = (\theta
_{1,0},\ldots,\theta_{d,0})'$ and $\bft_1 =
(\theta_{1,1},\ldots,\theta_{d,1})'$. It follows that the probability
that the logarithm of the Bayes factor exceeds a threshold
$\log(\gamma)$ can be expressed as
%
%
\begin{eqnarray}\label{dvec}
&&
\bP \bigl[ \log(\mathrm{BF}_{10}) > \log(\gamma) \bigr]\nonumber \\
&&\qquad= \bP \bigl\{ \bigl[
\eta_1(\theta_{1,1})-\eta_1(
\theta_{1,0}) \bigr] T_1({\mathbf x}) - \bigl[A(
\bft_1) - A(\bft_0) \bigr] > \log(\gamma) \bigr\}
\\
&&\qquad= \bP \biggl[ T_1({\mathbf x}) > \frac{\log(\gamma) +[A(\bft_1) -
A(\bft_0) ] }{ [\eta_1(\theta_{1,1})-\eta_1(\theta_{1,0})]}
\biggr].\nonumber
\end{eqnarray}

The probability in (\ref{dvec}) is maximized by minimizing the
right-hand side of the inequality. The extension to composite
alternative hypotheses follows the logic described in inequalities
(\ref{dvec1})--(\ref{dvec2}), which shows that UMPBT($\gamma$) tests
can be obtained in this setting by choosing the prior distribution of
$\bft_1$ under the alternative hypotheses so that it concentrates its
mass on the set
%
%
\begin{equation}
\argmin_{\theta}\frac{\log(\gamma) +[A(\bft_1) - A(\bft_0) ] }{
[\eta_1(\theta_{1,1})-\eta_1(\theta_{1,0})]},
\end{equation}
while maintaining the constraint that the values of $\eta_i(\bft)$
are equal under both hypotheses. Similar constructions apply if $\eta
_1$ is monotonically decreasing, or if the alternative hypothesis
specifies that $\theta_{1,0}<\theta_{0,0}$.

More practically useful extensions of UMPBTs can be obtained when it is
possible to integrate out nuisance parameters in order to obtain a
marginal density for the parameter of interest that falls within the
class of exponential family of models. An important example of this
type occurs in testing whether a regression coefficient in a linear
model is zero.

\subsection{\texorpdfstring{Test of linear regression coefficient, $\sigma^2$ known}
{Test of linear regression coefficient, sigma2 known}}\label{sec41}\label{linregknown}

Suppose that
%
%
\begin{equation}
\bfy\sim N\bigl(\bfX\bfb, \sigma^2 {\mathbf I}_n \bigr),
\end{equation}
where $\sigma^2$ is known, $\bfy$ is an $n\times1$ observation
vector, $\bfX$ an $n\times p$ design matrix of full column rank and
$\bfb=(\beta_1,\ldots,\beta_p)'$ denotes a $p\times1$ regression
parameter. The null hypothesis is defined as $H_0\dvtx  \beta_p =0$. For
concreteness, suppose that interest focuses on testing whether $\beta
_p>0$, and that under both the null and alternative hypotheses, the
prior density on the first $p-1$ components of $\bfb$ is a
multivariate normal distribution with mean vector ${\mathbf0}$ and
covariance matrix $\sigma^2 \bS$. Then the marginal density of $\bfy
$ under $H_0$ is
%
%
\begin{equation}
m_0(\bfy) = \bigl(2\pi\sigma^2\bigr)^{-n/2}
\llvert \bS\rrvert ^{-
{1}/{2}} \llvert \bF\rrvert ^{-{1}/{2}}\exp \biggl(-
\frac
{R}{2\sigma^2} \biggr),
\end{equation}
where
%
%
\begin{equation}
\label{RF}\quad \bF= \bfXp'\bfXp+ \bS^{-1},\qquad {\mathbf H} = \bfXp
\bF^{-1} \bfXp ',\qquad R = \bfy' ({\mathbf
I}_n-{\mathbf H})'\bfy,
\end{equation}
and $\bfXp$ is the matrix consisting of the first $p-1$ columns of
$\bfX$.\vadjust{\goodbreak}

Let $\beta_{p*}$ denote the value of $\beta_p$ under the alternative
hypothesis $H_1$ that defines the UMPBT($\gamma$), and let $\xp$
denote the $p$th column of $\bfX$. Then the marginal density of $\bfy
$ under $H_1$ is
%
%
\begin{equation}\qquad
m_1(\bfy) = m_0(\bfy) \times \exp \biggl\{ -
\frac{1}{2\sigma^2} \bigl[ \bps^2 \xp'({\mathbf
I}_n-{\mathbf H}) \xp- 2\bps\xp' ({\mathbf I}_n-{\mathbf
H}) \bfy \bigr] 
\biggr\}.
\end{equation}

It follows that the probability that the Bayes factor $\mathrm{BF}_{10}$ exceeds
$\gamma$ can be expressed as
%
%
\begin{equation}
\bP \biggl[ \xp'(\bfI_n -{\mathbf H})\bfy>
\frac{\sigma^2\log(\gamma
)}{\bps} + \frac{1}{2} \bps\xp'({\mathbf
I}_n-{\mathbf H}) \xp \biggr],
\end{equation}
which is maximized by minimizing the right-hand side of the inequality.
The UMPBT($\gamma$) is thus obtained by taking
%
%
\begin{equation}
\label{soln} \bps= \sqrt{\frac{2\sigma^2 \log(\gamma)}{\xp'({\mathbf I}_n-{\mathbf
H}) \xp}}.
\end{equation}
The corresponding one-sided test of $\beta_p<0$ is obtained by
reversing the sign of $\bps$ in (\ref{soln}).

Because this expression for the UMPBT assumes that $\sigma^2$ is
known, it is not of great practical significance by itself. However,
this result may guide the specification of alternative models in, for
example, model selection algorithms in which the priors on regression
coefficients are specified conditionally on the value of $\sigma^2$.
For example, the mode of the nonlocal priors described in \citet
{johnson12} might be set to the UMPBT values after determining an
appropriate value of $\gamma$ based on both the sample size $n$ and
number of potential covariates $p$.

\section{Approximations to UMPBTs using data-dependent alternatives}\label{sec5}
In some situations---most notably in linear models with unknown
variances---data dependent alternative hypotheses can be defined to
obtain tests that are approximately uniformly most powerful in
maximizing the probability that a Bayes factor exceeds a threshold.
This strategy is only attractive when the statistics used to define the
alternative hypothesis are ancillary to the parameter of interest.

\subsection{\texorpdfstring{Test of normal mean, $\sigma^2$ unknown}
{Test of normal mean, sigma2 unknown}}\label{sec51}

Suppose that $x_i$, $i=1,\ldots,n$, are i.i.d. $N(\mu,\sigma^2)$, that
$\sigma^2$ is unknown and that the null hypothesis is $H_0\dvtx  \mu= \mu
_0$. For convenience, assume further that the prior distribution on
$\sigma^2$ is an inverse gamma distribution with parameters $\alpha$
and $\lambda$ under both the null and alternative hypotheses.

To obtain an approximate UMPBT($\gamma$), first marginalize over
$\sigma^2$ in both models. Noting that
$(1+a/t)^t \rightarrow e^a $, it follows that the Bayes\vadjust{\goodbreak} factor in favor
of the alternative hypothesis satisfies
%
%
\begin{eqnarray}
\mathrm{BF}_{10}({\mathbf x}) &=& \biggl[ \frac{ \sum_{i=1}^n (x_i-\mu_0)^2 +
2\lambda}{ \sum_{i=1}^n (x_i-\mu_1)^2 + 2\lambda} \biggr]^{n/2+\alpha}
\\
&\approx& \biggl[ \frac{ 1 + {(\bar{{\mathbf x}} - \mu_0)^2}/{s^2}}{ 1 +
{(\bar{{\mathbf x}} - \mu_1)^2}/{s^2} } \biggr]^{n/2+\alpha}
\\
\label{tbf}
&\approx& \exp \biggl\{ -\frac{n}{2s^2} \bigl[ (\bar{{\mathbf x}}-\mu
_1)^2 - (\bar{{\mathbf x}}-\mu_0)^2
\bigr] \biggr\},
\end{eqnarray}
where
%
%
\begin{equation}
s^2 = \frac{ \sum_{i=1}^n (x_i-\bar{{\mathbf x}})^2 + 2\lambda
}{n+2\alpha}.
\end{equation}
The expression for the Bayes factor in (\ref{tbf}) reduces to (\ref
{normalBF}) if $\sigma^2$ is replaced by $s^2$. This implies that an
approximate, but data-dependent UMPBT alternative hypothesis can be
specified by taking
%
%
\begin{equation}
\label{adhoc} \mu_1 = \mu_0 \pm s \sqrt{
\frac{2\log\gamma}{n}},
\end{equation}
depending on whether $\mu_1>\mu_0$ or $\mu_1<\mu_0$.

\begin{figure}

\includegraphics{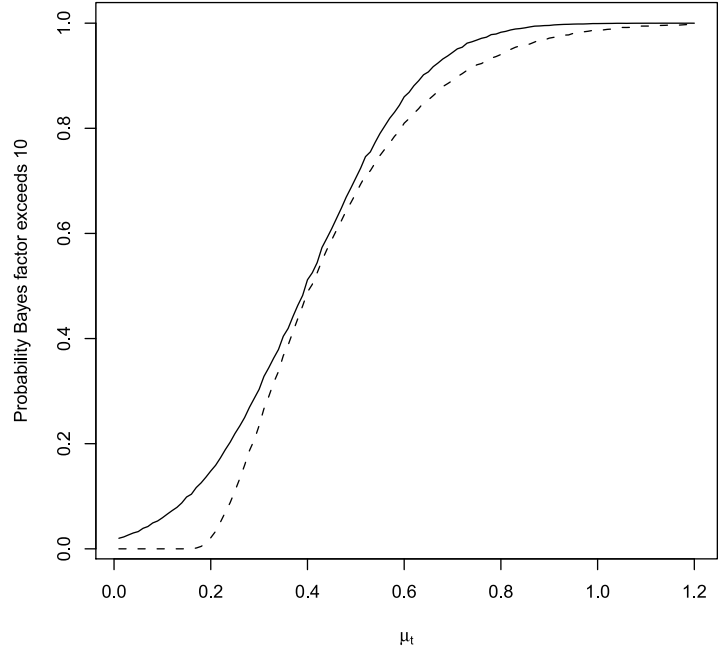}

\caption{Probability that Bayes factor based on data-dependent,
approximate UMPBT alternative exceeds 10 when $\mu_0=0$ and $n=30$
(solid curve). The dashed curve displays this probability when the
Bayes factor is calculated under the alternative hypothesis that $\mu
_1$ equals data-generating parameter (displayed on the horizontal axis)
and $\sigma^2=1$ (the true value).}
\label{unknownplot}
\end{figure}

Figure~\ref{unknownplot} depicts the probability that the Bayes factor
exceeds $\gamma=10$ when testing a null hypothesis that $\mu=0$ based
on an independent sample of size $n=30$ normal observations with unit
variance ($\sigma^2=1$) and using (\ref{adhoc}) to set the value of
$\mu_1$ under the alternative hypothesis. For comparison, the
probability that the Bayes factor exceeds 10 when the alternative is
defined by taking $\sigma^2=1$ and $\mu_1$ to be the data-generating
parameter is depicted by the dashed curve in the plot. Interestingly,
the data-dependent, approximate UMPBT(10) provides a higher probability
of producing a Bayes factor that exceeds 10 than do alternatives fixed
at the data generating parameters.

\subsection{\texorpdfstring{Test of linear regression coefficient, $\sigma^2$ unknown}
{Test of linear regression coefficient, sigma2 unknown}}\label{sec52}
As final example, suppose that the sampling model of Section~\ref{linregknown} holds, but assume now that the observational variance
$\sigma^2$ is unknown and assumed under both hypotheses to be drawn
from an inverse gamma distribution with parameters $\alpha$ and
$\lambda$. Also assume that the prior distribution for the first $p-1$
components of $\bfb$, given $\sigma^2$, is a multivariate normal
distribution with mean ${\mathbf0}$ and covariance matrix $\sigma^2
\Sigma$. As before, assume that $H_0\dvtx  \beta_p = 0$. Our goal is to
determine a value $\bps$ so that $H_1\dvtx  \beta_p = \bps$ is the
UMPBT($\gamma$) under the constraint that $\beta_p>0$.

Define $\bfy_1 = \bfy-\xp\bps$ and let $\bfy_0 = \bfy$. By
integrating with respect to the prior densities on $\sigma^2$ and the
first $p-1$ components of $\bfb$, the marginal density of the data
under hypothesis $i$, $i=0,1$ can be expressed as
%
%
\begin{equation}
m_i(\bfy) = 2^\alpha\pi^{-n/2} |
\Sigma|^{-1/2} \frac{\lambda
^\alpha}{\Gamma(\alpha)} \Gamma(n/2+\alpha) |{\mathbf
F}|^{-1/2} R_i^{-n/2-\alpha},
\end{equation}
where ${\mathbf F}$ is defined in (\ref{RF}), and
%
%
\begin{equation}
R_i = \bfy_i' ({\mathbf I}_n - {
\mathbf H}) \bfy_i + 2\lambda.
\end{equation}

It follows that the Bayes factor in favor of $H_1$ can be written as
%
%
\begin{eqnarray}
\mathrm{BF}_{10} &=& \biggl[ 1+ \frac{\bps^2 \xp' ({\mathbf I}_n-{\mathbf H}) \xp
- 2\bps\xp' ({\mathbf I}_n-{\mathbf H}) \bfy}{R_0} \biggr]^{-n/2-\alpha
}
\\
\label{ee}
&\approx& \exp \biggl\{ -\frac{1}{2s_p^2} \bigl[ \bps^2
\xp' ({\mathbf I}_n-{\mathbf H}) \xp- 2\bps\xp' ({
\mathbf I}_n-{\mathbf H}) \bfy \bigr] \biggr\},
\end{eqnarray}
where
%
%
\begin{equation}
s_p^2 = \frac{R_0 }{n+2\alpha}.
\end{equation}

The UMPBT($\gamma$) is defined from (\ref{ee}) according to
%
%
\begin{equation}\quad
\bP (\mathrm{BF}_{10} > \gamma ) = \bP \biggl[ \xp'({\mathbf
I}_n-{\mathbf H}) \bfy> \frac{s_p^2 \log(\gamma)}{\bps} + \frac{1}{2} \bps\xp
' ({\mathbf I}_n-{\mathbf H}) \xp \biggr].
\end{equation}
Minimizing the right-hand side of the last inequality with respect to
$\bps$ results in
%
%
\begin{equation}
\bps= \sqrt{ \frac{2s_p^2 \log(\gamma)}{\xp' ({\mathbf I}_n-{\mathbf H})
\xp}}.
\end{equation}

This expression is consistent with the result obtained in the known
variance case, but with $s_p^2$ substituted for $\sigma^2$.

\section{Discussion}\label{sec6}\label{discussion}

The major contributions of this paper are the definition of UMPBTs and
the explicit description of UMPBTs for regular one-parameter
exponential family models. The existence of UMPBTs for exponential
family models is important because these tests represent the most
common hypothesis tests conducted by practitioners. The availability of
UMPBTs for these models means that these tests can be used to interpret
test results in terms of Bayes factors and posterior model
probabilities in a wide range of scientific settings. The utility of
these tests is further enhanced by the connection between UMPBTs and
UMPTs that have the same rejection region. This connection makes it
trivial to simultaneously report both the $p$-value from a test and the
corresponding Bayes factor.

The simultaneous report of default Bayes factors and $p$-values may
play a pivotal role in dispelling the perception held by many scientists
that a $p$-value of 0.05 corresponds to ``significant'' evidence
against the null hypothesis.
The preceding sections contain examples in which this level of
significance favors the alternative hypothesis by odds of only 3 or 4
to 1.
Because few researchers would regard such odds as strong evidence in
favor of a new theory, the use of UMPBTs and the report of Bayes
factors based upon them may lead to more realistic interpretations of
evidence obtained from scientific studies.

The large sample properties of UMPBTs described in Section~\ref{sec21} deserve
further comment. From Lemma~\ref{lem2}, it follows that the expected weight of
evidence in favor of a true null hypothesis in an exponential family
model converges to $\log(\gamma)$ as the sample size $n$ tends to
infinity. In other words, the evidence threshold $\gamma$ represents
an approximate bound on the evidence that can be collected in favor of
the null hypothesis.
This implies that $\gamma$ must be increased with $n$ in order to
obtain a consistent sequence of tests.

Several criteria might be used for selecting a value for $\gamma$ in
large sample settings. One criterion can be inferred from the first
statement of Lemma~\ref{lem2}, where it is shown that the difference between the
tested parameter's value under the null and alternative hypotheses is
proportional to $[\log(\gamma)/n]^{1/2}$. For this difference to be a
constant---as it would be in a subjective Bayesian test---$\log(\gamma
)$ must be proportional to $n$, or $\gamma= \exp(cn)$ for some
$c>0$. This suggests that an appropriate value for $c$ might be
determined by calibrating the weight of evidence against an accepted
threshold/sample size combination. For example, if an evidence
threshold of 4 were accepted as the standard threshold for tests
conducted with a sample size of 100, then $c$ might be set to $\log
(4)/100 = 0.0139$. This value of $c$ leads to an evidence threshold of
$\gamma= 16$ for sample sizes of 200, a threshold of 64 for sample
sizes of 300, etc. From (\ref{gammaequiv}), the significance levels
for corresponding $z$-tests would be 5\%, 1\% and 0.2\%, respectively.

The requirement to increase $\gamma$ to achieve consistent tests in
large samples also provides insight into the performance of standard
frequentist and subjective Bayesian tests in large sample settings. The
exponential growth rate of $\gamma$ required to maintain a fixed
alternative hypothesis suggests that the weight of evidence should be
considered against the backdrop of sample size, even in Bayesian tests.
This is particularly important in goodness-of-fit testing where small
deviations from a model may be tolerable. In such settings, even
moderately large Bayes factors against the null hypotheses may not be
scientifically important when they are based on very large sample sizes.

From a frequentist perspective, the use of UMPBTs in large sample
settings can provide insight into the deviations from null hypotheses
when they are (inevitably) detected. For instance, suppose that a
one-sided 1\% test has been conducted to determine if the mean of
normal data is 0, and that the test is rejected with a $p$-value of
0.001 based on a sample size of 10,000. From (\ref{gammaequiv}), the
implied evidence threshold for the test is $\gamma=15$, and the
alternative hypothesis that has been implicitly tested with the UMPBT
is that $\mu=0.023\sigma$. Based on the observation of $\bar
{x}=0.031\sigma$, the Bayes factor in favor of this alternative is
88.5. Although there are strong odds against the null, the scientific
importance of this outcome may be tempered by the fact that the
alternative hypothesis that was supported against the null represents a
standardized effect size of only 2.3\%.

This article has focused on the specification of UMPBTs for one-sided
alternatives. A simple extension of these tests to two-sided
alternatives can be obtained by assuming that the alternative
hypothesis is represented by two equally-weighted point masses located
at the UMPBT values determined for one-sided tests. The Bayes factors
for such tests can be written as
%
%
\begin{equation}
\bP \biggl[\frac{0.5 m_{l}({\mathbf x})+0.5 m_{h}({\mathbf x})}{m_0({\mathbf x})} >\gamma \biggr],
\end{equation}
where $m_{l}$ and $m_{h}$ denote marginal densities corresponding to
one-sided UMPBTs. Letting $m_{*}({\mathbf x}) = \max(m_{l}({\mathbf
x}),m_{h}({\mathbf x}))$ for the data actually observed, and assuming that
the favored marginal density dominates the other, it follows that
%
%
\begin{equation}
\bP \biggl[\frac{0.5 m_{l}({\mathbf x})+0.5 m_{h}({\mathbf x})}{m_0({\mathbf x})
}>\gamma \biggr] \approx \bP \biggl[
\frac{ m_{*}({\mathbf x})}{m_0({\mathbf x}) }>2\gamma \biggr].
\end{equation}
Thus, an approximate two-sided UMPBT($\gamma$) can be defined by
specifying an alternative hypothesis that equally concentrates its mass
on the two one-sided UMPBT($2\gamma$) tests.

Additional research is needed to identify classes of models and testing
contexts for which UMPBTs can be defined. The UMPBTs described in this
article primarily involve tests of point null hypotheses, or tests that
can be reduced to a test of a point null hypothesis after marginalizing
over nuisance parameters. Whether UMPBTs can be defined in more general
settings remains an open question.

\section*{Acknowledgments}

The author thanks an Associate Editor and two referees for numerous
comments that improved this article. Article
content is solely the responsibility of the author and does not
necessarily represent the official views of the National Cancer
Institute or the National Institutes of Health.



\printaddresses

\end{document}